\documentstyle[amssymb,12pt]{article}

\def\oisi{\buildrel\textstyle{\rm i \atop ^\vee} \over.}
\def\oisk{\buildrel\textstyle{\rm k \atop ^\vee} \over\dots}
\def\oisr{\buildrel\textstyle{\rm r \atop ^\vee} \over\dots}
\def\oiss{\buildrel\textstyle{\rm s \atop ^\vee} \over\dots}
\def\Talpha#1{\vbox{\ialign{##\crcr
    $\alpha$\crcr\noalign{\kern2pt\nointerlineskip}
           $\hfil\displaystyle{#1}\hfil$\crcr}}}

\def\Onabla#1{\vbox{\ialign{##\crcr
   $\,\scriptstyle{0}$\crcr\noalign{\kern2pt\nointerlineskip}
           $\hfil\displaystyle{#1}\hfil$\crcr}}}

\def\cala{{\cal A}}

\def\calb{{\cal B}}
\def\calc{{\cal C}}

\def\calr{{\cal R}}

\def\bbbone{\mbox{\rm 1\hspace {-.6em} l}}
\def\ham{{\mbox{Ham}}}
\def\hom{{\mbox{Hom}}}

\def\der{{\mbox{\scriptsize Der}}}
\def\gder{{\mbox{Der}}}
\def\os{\underline{\Omega}}

\def\fc{{\frak C}}

\begin{document}

\baselineskip=0.7cm
\begin{center}
{\Large\bf Some Aspects of Noncommutative Differential Geometry}
\end{center}
\vspace{0.75cm}

\begin{center}
Michel DUBOIS-VIOLETTE

\vspace{0.3cm}
{\small
Laboratoire de Physique Th\'eorique et Hautes
Energies\footnote{Laboratoire associ\'e au Centre National de la
Recherche Scientifique - URA D0063}\\ Universit\'e Paris XI, B\^atiment
211\\ 91 405 Orsay Cedex, France\\ flad@qcd.th.u-psud.fr}

\end{center} \vspace{1cm}

\begin{center} November 1995 \end{center}

\vspace {5cm}

\noindent L.P.T.H.E.-ORSAY 95/78\\
\noindent ESI-preprint

\newpage \begin{abstract}
We discuss in some generality aspects of noncommutative differential geometry
associated with reality conditions and with differential calculi. We then
describe the differential calculus based on derivations as generalization of
vector fields, and we show its relations with quantum mechanics. Finally we
formulate a general theory of connections in this framework.
\end{abstract}

\section{Introduction}

In \cite{koszul:cal}, J.L. Koszul described a powerful algebraic version of
differential geometry in terms of a commutative associative algebra $\calc$,
$\calc$-modules and connections (``derivation laws") on these modules. For the
applications to differential geometry, $\calc$ is the algebra of smooth
functions on a manifold and the $\calc$-modules are modules of smooth sections
of smooth vector bundles over the manifold. The fact that classical
differential geometry admits such an algebraic formulation is at the very
origin of the idea of noncommutative differential geometry. Historically, the
motivation of noncommutative geometry was the development of quantum theory
\cite{pamd}. In noncommutative geometry, one replaces the commutative
associative algebra $\calc$ by an associative algebra $\cala$ which is not
assumed to be commutative. However this replacement raises several problems
which will be discussed in this lecture.\\

First problem: what should replace the $\calc$-modules? The problem arises
because there are at least four inequivalent generalizations of the notion of a
module over a commutative algebra when the algebra is replaced by a
noncommutative algebra $\cala$. There is the notion of right $\cala$-module and
the dual notion of left $\cala$-module. If one recalls that a module over a
commutative algebra is canonically a bimodule (of a specific kind), there is a
notion of bimodule over $\cala$. Finally, since a commutative algebra coincides
with its center, there is the notion of module over the center $Z(\cala)$ of
$\cala$. As will be explained latter, there is also a duality between
$Z(\cala)$-modules and bimodules over $\cala$.\\

Second problem: what should be the generalization of the classical notions of
reality? For classical differential geometry one can use for $\calc$ either the
real commutative algebra of smooth real-valued functions or the complex
commutative $\ast$-algebra of smooth complex-valued  functions. More generally,
if $\calc$ is a complex commutative $\ast$-algebra then the set $\calc^h$ of
its hermitian elements is a real commutative algebra and $\calc$ is the
complexification of $\calc^h$. Conversely if $\calc_{\Bbb R}$ is a real
commutative algebra, then its complexification $\calc$ is canonically a complex
commutative $\ast$-algebra and one has $\calc_{\Bbb R}=\calc^h$. In fact
$\calc\mapsto \calc^h$ is an equivalence of the category of commutative
associative $\ast$-algebras over $\Bbb C$ and $\ast$-homomorphisms onto the
category of commutative associative algebras over $\Bbb R$ and homomorphisms of
real algebras. The situation is quite different for noncommutative algebras. If
$\cala$ is a complex associative $\ast$-algebra, the set $\cala^h$ of its
hermitian elements is generally not an associative algebra but a real Jordan
algebra. This means that one has two choices for the generalization of the
algebra of real-valued functions, either the real Jordan algebra $\cala^h$ of
all hermitian element of a complex associative $\ast$-algebra $\cala$, which
plays the role of the algebra of complex-valued functions, or a real
associative algebra. Here we take the first point of view. This choice, which
is the standard one, is dictated by quantum theory and, more generally, by
spectral theory. This reality problem is not independent of the first problem
because if $\calc$ is a complex commutative associative $\ast$-algebra there is
again an obvious equivalence between the involutive $\calc$-modules and the
$\calc^h$-modules \cite{mdv:pm1}.\\

Third problem: which differential calculus should be used? In other words what
should be the generalization of differential forms? Such a generalization is
needed for instance to define connections. There is a minimal set of
assumptions which must be satisfied and which will be described in the sequel
but nevertheless the choice is not straightforward. We can make here the
following remarks. In his pioneer work on the subject \cite{connes:02}, A.
Connes defined the cyclic cohomology of an algebra and showed that the correct
generalization of the homology of a manifold is the reduced cyclic cohomology.
This means that the generalization of the cohomology of a manifold in
noncommutative geometry must be the reduced cyclic homology of the algebra
$\cala$ which replaces the algebra of smooth functions. In classical
differential geometry, the de Rham theorem states that the cohomology of a
manifold, (a topological invariant), coincides with the cohomology of its
differential forms. This does not mean that any cochain complex which has the
reduced cyclic homology as cohomology is an acceptable generalization of
differential forms, and this for at least two reasons. First, even in the
classical situation, there are many ways to compute the cohomology of a
manifold and, in particular, there are complexes which are not connected with
the differential structure and which have this cohomology. Second, the de Rham
theorem is not a tautological result but a deep theorem of differential
topology which means that there may well be proper noncommutative
generalizations of differential geometry for which the generalization of de
Rham theorem fails to be true.\\

The problems quoted above will be discussed in the first part of this lecture.
Then the differential calculus based on the derivations as generalization of
vector fields will be introduced \cite{dv:1}. This differential calculus is the
direct generalization of the one used by J.L. Koszul in \cite{koszul:cal}; it
is also connected with the differential calculus used by A.~Connes for
noncommutative dynamical systems in \cite{connes:03}. The noncommutative
symplectic structures will be defined in this framework and the relation with
quantum mechanics will be described. Finally we shall describe  the theory of
connections in this framework. Examples of such connections and applications to
gauge field theory may be found in \cite{dv:2}, \cite{dvkm:1}, \cite{dvkm:2},
\cite{mdv:pm1}, \cite{mad}. \\

Let $\cala$ be an associative algebra. If $M$  and $N$ are right
$\cala$-modules, the space of all right $\cala$-module homomorphisms of $M$
into $N$ will be denoted by $\hom^{\cala}(M,N)$; if $M$ and $N$ are left
$\cala$-modules, the space of all left $\cala$-module homomorphisms of $M$ into
$N$ will be denoted by $\hom_{\cala}(M,N)$. When $\cala$ is a commutative
algebra $\calc$, both notions coincide and the space of $\calc$-module
homomorphisms of $M$ into $N$ will be denoted by $\hom_{\calc}(M,N)$. If
$\calb$ is another associative algebra and if $M$ and $N$ are
$(\cala,\calb)$-bimodules,  $\hom^{\calb}_{\cala}(M,N)$ will denote the space
of all bimodules homomorphisms of $M$ into $N$. In the sequel, we shall often
use the word algebra to mean associative algebra.\\

This lecture is partly based on joint works with R. Kerner, J. Madore and P.W.
Michor
%% FOLLOWING LINE CANNOT BE BROKEN BEFORE 80 CHAR
%% FOLLOWING LINE CANNOT BE BROKEN BEFORE 80 CHAR
\cite{dv:1},\cite{dv:2},\cite{dvkm:1},\cite{dvkm:2},\cite{mdv:pm1},\cite{mdv:pm2},\cite{mdv:pm3} and the author is grateful to John Madore for discussions and careful reading of the manuscript.

\section{Modules, bimodules and reality}

In the following $\cala$ is a complex unital associative $\ast$-algebra which
is to be considered as a noncommutative generalization of an algebra of complex
functions. As a consequence what must replace the algebra of real-valued
functions is generally not an associative real algebra but the Jordan algebra
$\cala^h$ of hermitian elements of $\cala$. Although quite familiar in quantum
theory, this fact has non trivial consequences for the noncommutative
generalization of classical reality conditions. In fact we are here interested
in noncommutative differential geometry, which means that $\cala$ is to be
considered as the generalization of the algebra of complex smooth functions on
a manifold. The Jordan algebra $\cala^h$ replaces then the algebra of real
smooth functions.\\

Let $E$ be a smooth complex vector bundle of finite rank over a manifold $V$.
Then the set $\Gamma(E)$ of its smooth sections is a finite projective module
over the algebra $\calc^\infty(V)$ of smooth complex functions on $V$.
Furthermore the correspondence $E\mapsto\Gamma(E)$ is an equivalence of the
category of smooth complex vector bundles of finite rank over $V$ onto the
category of finite projective modules over $\calc^\infty(V)$. Let now $E_{\Bbb
R}$ be a smooth real vector bundle over $V$, its complexification $E$ is a
smooth complex vector bundle over $V$ equipped with a canonical antilinear
involution $\xi\mapsto\xi^\ast$ such that $\xi\in E_{\Bbb R}$ if and only if
$\xi=\xi^\ast$. The module $\Gamma(E)$ is then a $\ast$-module over the
$\ast$-algebra $\calc^\infty(V)$ in the sense that it is equipped with an
antilinear involution $\psi\mapsto\psi^\ast$ such that
$(f\psi)^\ast=f^\ast\psi^\ast,\ \ \forall f\in \calc^\infty(V)$ and $\forall
\psi\in \Gamma(E)$, where $f\mapsto f^\ast$ is the complex conjugation. A
section of $E_{\Bbb R}$ is a section $\psi\in \Gamma(E)$ such that
$\psi=\psi^\ast$. Clearly, one can replace $E_{\Bbb R}$ by the $\ast$-module
$\Gamma(E)$. With this in mind, let us more generally consider the notion of
module and the notion of $\ast$-module over a commutative $\ast$-algebra
$\calc$ and investigate their generalizations when $\calc$ is replaced by the
noncommutative $\ast$-algebra $\cala$.\\

As pointed out in the introduction a $\calc$-module has several natural
generalizations: a right $\cala$-module, a left $\cala$-module, a module over
the center $Z(\cala)$ of $\cala$ and a bimodule over $\cala$. Right
$\cala$-modules and left $\cala$-modules are dual in the sense that if $M$ is a
right $\cala$-module, its dual $M^\ast=\hom^{\cala}(M,\cala)$ is a left
$\cala$-module and if $N$ is a left $\cala$-module, its dual
$N^\ast=\hom_{\cala}(N,\cala)$ is a right $\cala$-module; this duality
generalizes the duality of $\calc$-modules. Similarily, there is a natural
duality between bimodules over $\cala$ and $Z(\cala)$-modules \cite{mdv:pm1}:
if $M$ is a bimodule over $\cala$, {\sl its $\cala$-dual}
$M^{\ast_{\cala}}=\hom^{\cala }_{\cala}(M,\cala)$ is canonically a
$Z(\cala)$-module and if $N$ is a $Z(\cala)$-module, {\sl its $\cala$-dual}
$N^{\ast_{\cala}}=\hom_{Z(\cala)}(N,\cala)$ is canonically a bimodule over
$\cala$. This duality ($\cala$-duality) also generalizes the duality of
$\calc$-modules when the bimodules over $\calc$ are the underlying bimodules of
$\calc$-modules.\\

Concerning the generalization of $\ast$-modules over $\calc$, (i.e. the
generalization of the description of real vector bundles), one notices that one
cannot use right or left $\cala$-modules because, since the involution of
$\cala$ reverses the order of the product in $\cala$, there cannot be a notion
of right or left $\ast$-module over $\cala$. In contrast, since $Z(\cala)$ is a
commutative $\ast$-algebra, the notion of $\ast$-module over $Z(\cala)$ is
perfectly defined and one can introduce a dual notion of $\ast$-bimodule over
$\cala$: a bimodule $M$ over $\cala$ is {\sl a $\ast$-bimodule} over $\cala$ if
it is equipped with an antilinear involution $m\mapsto m^\ast$ such that one
has $(xmy)^\ast=y^\ast m^\ast x^\ast\ \ \forall x,y\in \cala$ and $\forall m\in
M$.\\

Thus, simple considerations of reality rule out right or left $\cala$-modules
for the description of a generalization of real vector bundles. This does not
mean that one cannot use them for the generalization of complex vector bundles,
this simply means that all the above generalizations of the notion of
$\calc$-module have to be considered when $\calc$ is replaced by the
noncommutative algebra $\cala$. The fact that bimodule structures arise in
connection with reality in noncommutative geometry has been also pointed out in
\cite{connes:05} by A.~Connes in the context of his spectral triples approach
to noncommutative geometry \cite{connes:03},\cite{connes:04}.\\

It must be stressed that not every bimodule over a commutative algebra $\calc$
is the underlying bimodule of a $\calc$-module and therefore not every bimodule
over $\cala$ can be considered as the generalization of a $\calc$-module. One
must select an appropriate class of bimodules, for instance the class of
central bimodules \cite{mdv:pm1},\cite{mdv:pm2}. A bimodule $M$ over $\cala$ is
called a {\sl central bimodule} if one has $zm=mz,\ \ \forall m\in M$ and
$\forall z\in Z(\cala)$. A central bimodule over a commutative algebra $\calc$
is just a $\calc$-module for its underlying bimodule structure. In
\cite{mdv:pm2}, the more restrictive notion of diagonal bimodule was
introduced. A bimodule $M$ over $\cala$ is called a {\sl diagonal bimodule} if
it is isomorphic to a subbimodule of $\cala^I$ for some set $I$. A diagonal
bimodule is central. A bimodule $M$ over $\cala$ is diagonal if and only if the
canonical mapping of $M$ into its $\cala$-bidual $M^{\ast{_\cala}\ast{_\cala}}$
is injective. In particular a diagonal bimodule over a commutative algebra
$\calc$ is just a $\calc$-module such that the canonical mapping in its bidual
is injective; projective $\calc$-modules are therefore diagonal bimodules. If
$N$ is a $Z(\cala)$-module, its $\cala$-dual $N^{\ast_{\cala}}$ is a diagonal
bimodule over $\cala$.

\section{Differential calculus}

In this section we wish to discuss some general features of the noncommutative
versions of differential forms.\\

A {\sl graded differential $\ast$-algebra} is a complex graded differential
algebra $\Omega=\oplus_{n\in \Bbb N} \Omega^n$ equipped with an antilinear
involution $\omega\mapsto \omega^\ast$ which preserves the degree and satisfies
$(\alpha\beta)^\ast=(-1)^{ab}\beta^\ast\alpha^\ast$ and
$(d\omega)^\ast=d(\omega^\ast)$ for $\alpha\in\Omega^a, \beta\in\Omega^b$ and
$\omega\in \Omega$, where $d$ is the differential of $\Omega$. Notice that then
$\Omega^0$ is a $\ast$-algebra. Given (as before) the complex unital
$\ast$-algebra $\cala$, {\sl a differential calculus over} $\cala$ is a graded
differential $\ast$-algebra $\Omega$ with $\Omega^0=\cala$. Among the
differential calculi over $\cala$, there is a universal one \cite{kar},
$\Omega_u(\cala)$, which we now review.\\

Let $\mu:\cala\otimes\cala\rightarrow \cala$ be the product $\mu(x\otimes
y)=xy$. The mapping $\mu$ is a bimodule homomorphism so its kernel
$\Omega^1_u(\cala)$ is a bimodule over $\cala$. One defines a derivation $d_u$
of $\cala$ into $\Omega^1_u(\cala)$ by setting $d_ux=\bbbone \otimes x -
x\otimes \bbbone$ for $x\in \cala$. The pair $(\Omega^1_u(\cala),d_u)$ is
characterized uniquely (up to an isomorphism) by the following universal
property \cite{cart:02}, \cite{bour}: given a derivation
$\delta:\cala\rightarrow M$ of $\cala$ into a bimodule $M$ over $\cala$, there
is a unique bimodule homomorphism $j_\delta:\Omega^1_u(\cala)\rightarrow M$
such that $\delta=j_\delta\circ d_u$. Let $\Omega_u(\cala)$ be the tensor
algebra over $\cala$ of the bimodule $\Omega^1_u(\cala)$ i.e.
$\Omega^0_u(\cala)=\cala$ and
$\Omega^n_u(\cala)=\otimes^n_{\cala}\Omega^1_u(\cala)$ for $n\geq 1$. The
derivation $d_u$ extends uniquely into a differential, again denoted by $d_u$,
of the graded algebra $\Omega_u(\cala)$. Using the above universal property of
$(\Omega^1_u(\cala),d_u)$ and the universal property of the tensor product over
$\cala$, one sees that the graded differential algebra $\Omega_u(\cala)$ is
characterized by the following universal property: given a graded differential
algebra $\Omega=\otimes_n\Omega^n$ with $\Omega^0=\cala$, there is a unique
homomorphism of graded differential algebra $\varphi:\Omega_u(\cala)\rightarrow
\Omega$ which induces the identity mapping of $\cala$ onto itself (i.e.
$\varphi\restriction\cala=id_{\cala}$). Furthermore, there is a unique
antilinear involution $\omega\mapsto\omega^\ast$ on $\Omega_u(\cala)$ which
extends the involution of $\cala$ and for which it is a graded differential
$\ast$-algebra \cite{slw}; this involution is induced on
$\Omega^n_u(\cala)(\subset \otimes^{n+1} \cala)$ by the involution of
$\otimes^{n+1}\cala$ defined by $(x_0\otimes x_1\otimes \dots \otimes
x_n)^\ast=(-1)^{\frac{n(n+1)}{2}} x^\ast_n \otimes \dots \otimes
x^\ast_1\otimes x^\ast_0$. Equipped with this involution, $\Omega_u(\cala)$ is
a differential calculus over $\cala$ which is universal in the sense that for
any differential calculus $\Omega$ over $\cala$ there is a unique homomorphism
of graded differential $\ast$-algebra of $\Omega_u(\cala)$ into $\Omega$ which
induces the identity mapping of $\cala$ onto itself.\\

One can expect, and it is our point of view here, that a noncommutative
generalization of differential forms is a differential calculus over $\cala$
when $\cala$ replaces the algebra $\calc^\infty(V)$ of smooth functions on a
manifold $V$. However not every differential calculus over $\cala$ is
appropriate. For instance the universal differential calculus is not a proper
generalization of the algebra of differential forms. Indeed
$\Omega_u(\calc^\infty(V))$ does not coincide with the algebra $\Omega(V)$ of
differential forms on $V$ although, by the universal property, there is a
homomorphism of graded differential algebra of $\Omega_u(\calc^\infty(V))$ into
$\Omega(V)$. More generally, if $\calc$ is a commutative algebra, the bimodule
$\Omega^1_u(\calc)$, for instance, is not the underlying bimodule of a module
since left and right multiplications by elements of $\calc$ do not coincide. In
any case the choice of a differential calculus $\Omega$ over $\cala$ as
generalization of the algebra of complex differential forms is not unique and
depends on the applications one has in mind \cite{connes:02}, \cite{connes:03},
\cite{rcoq}, \cite{dv:1}, \cite{mdv:pm1}, \cite{mdv:pm3}, \cite{kast},
\cite{mad},  \cite{slw}. In the next section we will describe a choice for
$\Omega$ based on derivations as generalization of vector fields. This choice,
which is a direct generalization of \cite{koszul:cal}, is natural in the sense
that it only depends on the algebra $\cala$ (and not on additional
structures).\\

Before leaving this section, two points are worth noticing. First some authors,
e.g. G.~Maltsiniotis \cite{malt}, consider that a proper generalization of
differential geometry is given by a graded differential algebra which then
replaces the algebra of differential forms; this point of view is more general
than the one, implicit here, where $\cala$ replaces the algebra of smooth
functions. Second there are generalizations of the space of differential forms
which are not differential algebras but merely differential complexes. For
instance, it was shown in \cite{kar} that the subspace
$[\Omega_u(\cala),\Omega_u(\cala)]$ of graded commutators in $\Omega_u(\cala)$
is stable by $d_u$ and that the cohomology of the complex
$(\Omega_u(\cala)/[\Omega_u(\cala),\Omega_u(\cala)],d_u)$, (which is not a
differential algebra in general), is the reduced cyclic homology of $\cala$
which in many aspects is a good generalization of de Rham cohomology. This is
why this complex is a natural generalization of the de Rham complex which is
often called {\sl the noncommutative de Rham complex}.

\section{Derivations and differential calculus}

In this section we explain our approach to the (noncommutative) differential
calculus over $\cala$, (a complex unital $\ast$-algebra), based on the
derivations of $\cala$ as generalization of vector fields \cite{dv:1},
\cite{dv:2}, \cite{dvkm:1}, \cite{dvkm:2}, \cite{mdv:pm1}, \cite{mdv:pm2},
\cite{mdv:pm3}. This approach is a noncommutative generalization of the one of
J.L. Koszul \cite{koszul:cal} which is based	 on the fact that a vector field
on a manifold $V$, i.e. a smooth section of the tangent bundle over $V$, is the
same thing as a derivation of the algebra $\calc^\infty(V)$ of smooth functions
on $V$. More generally, since the derivations are the infinitesimal algebra
automorphisms, they are the natural right-hand sides of differential evolution
equations. This is why the differential calculus based on derivations is the
natural one for commutative and noncommutative dynamical systems i.e. for
classsical as well as for quantum mechanics.\\

Let $\gder(\cala)$ denote the space of all derivations of $\cala$, i.e. the
space of all linear mappings $X$ of $\cala$ into itself satisfying the Leibniz
rule $X(xy)=X(x)y+xX(y)$. The space $\gder(\cala)$ is in a natural way a module
over the center $Z(\cala)$ of $\cala$ and in fact a $\ast$-module over
$Z(\cala)$ when equipped with the involution $X\mapsto X^\ast$ defined by
$X^\ast(x)=(X(x^\ast))^\ast$. The space $\gder(\cala)$ is also a Lie algebra
with Lie bracket $(X,Y)\mapsto [X,Y]=X\circ Y-Y\circ X$. This bracket satisfies
the reality condition $[X,Y]^\ast=[X^\ast,Y^\ast]$. Furthermore, $Z(\cala)$ is
stable under $\gder(\cala)$ and one has $[X,zY]=X(z)Y+z[X,Y]$, for any $X,Y\in
\gder(\cala)$ and $z\in Z(\cala)$. This last equality ensures that, in the
complex $C(\gder(\cala),\cala)$ of the $\cala$-valued Lie-algebra cochains of
$\gder(\cala)$, the subspace $\os_{\der}(\cala)$ of $Z(\cala)$-multilinear
cochains is stable under the differential, i.e. is a subcomplex.\\

More precisely, let $\os^n_{\der}(\cala)$ be the space of
$Z(\cala)$-multilinear antisymmetric mappings of $(\gder(\cala))^n$ into
$\cala$, (i.e. $\os^n_{\der}(\cala)=\hom_{Z(\cala)}(\Lambda^n_{Z(\cala)}
\gder(\cala),\cala))$. Then the graded space $\os_{\der}(\cala)=\oplus_n
\os^n_{\der}(\cala)$ is in a natural way a graded algebra (the product
combining the product of $\cala$ with antisymmetrisation in the arguments). One
verifies that one defines a differential $d$ of $\os_{\der}(\cala)$, i.e. an
antiderivation of degree 1 satisfying $d^2=0$, by setting, for
$\omega\in\os^n_{\der}(\cala)$ and $X_i\in\gder (\cala)$,
$$(d\omega)(X_0,\dots,X_n)=\sum^n_{k=0}(-1)^k X_k \omega(X_0,\oisk,X_n)$$
$$+\sum_{0\leq r<s\leq n}(-1)^{r+s} \omega([X_r,X_s],X_0,\oisr\oiss,X_n)$$
where $\oisi$ means omission of $X_i$. Thus, equipped with this differential,
$\os_{\der}(\cala)$ is a graded differential algebra and the subalgebra
$\os^0_{\der}(\cala)$ coincides with $\cala$. If one equips $\os_{\der}(\cala)$
with the involution $\omega\mapsto\omega^\ast$ defined by $\
\omega^\ast(X_1,\dots,X_n)\ =\\
(\omega(X^\ast_1,\dots,X^\ast_n))^\ast,$ it becomes a differential calculus
over $\cala$.\\

Let $\Omega_{\der}(\cala)$ be the smallest differential subalgebra of
$\os_{\der}(\cala)$ which contains $\cala$. The differential algebra
$\Omega_{\der}(\cala)$ is the canonical image of $\Omega_u(\cala)$ in
$\os_{\der}(\cala)$ and is stable by the involution; it consists of finite sums
of elements of the form $x_0dx_1\dots dx_n$, $x_i\in\cala$. The graded
differential $\ast$-algebra $\Omega_{\der}(\cala)$ is also a differential
calculus over $\cala$.\\

Both $\Omega_{\der}(\cala)$ and $\os_{\der}(\cala)$ are generalizations of the
algebra of complex differential forms. If $V$ is a finite-dimensional
paracompact manifold then $\Omega_{\der}(\calc^\infty(V))$ and
$\os_{\der}(\calc^\infty(V))$ both coincide with the graded differential
$\ast$-algebra $\Omega(V)$ of complex differential forms on $V$. In general the
inclusion $\Omega_{\der}(\cala)\subset \os_{\der}(\cala)$ is a strict one:
$\Omega_{\der}(\cala)$ is the minimal version of noncommutative differential
forms based on derivations while $\os_{\der}(\cala)$ is the maximal one. It is
worth noticing here that even in the classical situation the above inclusion
may be strict, e.g. if $V$ is a manifold which does not admit a partition of
unity then the inclusion $\Omega_{\der}(\calc^\infty(V))\subset
\os_{\der}(\calc^\infty(V))$ is a strict one. There is however a density result
of $\Omega_{\der}(\cala)$ in $\os_{\der}(\cala)$ which we now describe at the
level of one-forms \cite{mdv:pm1}.\\

By its very definition, the bimodule $\os^1_{\der}(\cala)$ is the $\cala$-dual
of the $Z(A)$-module $\gder(\cala)$, i.e.
%% FOLLOWING LINE CANNOT BE BROKEN BEFORE 80 CHAR
%% FOLLOWING LINE CANNOT BE BROKEN BEFORE 80 CHAR
$\os^1_{\der}(\cala)=(\gder(\cala))^{\ast_{\cala}}=\hom_{Z(\cala)}(\gder(\cala),\cala)$. On the other hand, by the universal property of $(\Omega^1_u(\cala),d_u)$, $\gder(\cala)$ can be identified with $\hom^{\cala}_{\cala}(\Omega^1_u(\cala),\cala)$ through the canonical mapping $X\mapsto j_X$ (see in last section). However the intersection of the kernels of the bimodule homomorphisms of $\Omega^1_u(\cala)$ into $\cala$, (which is the intersection of the kernels of the $j_X$ when $X$ runs over $\gder(\cala)$), is just the kernel of the canonical bimodule homomorphism of $\Omega^1_u(\cala)$ onto $\Omega^1_{\der}(\cala)$ \cite{dv:1} and therefore one has $\hom^{\cala}_{\cala}(\Omega^1_u(\cala),\cala)=
\hom^{\cala}_{\cala}(\Omega^1_{\der}(\cala),\cala)$. So one has finally
$\hom^{\cala}_{\cala}(\Omega^1_{\der}(\cala),\cala)=\gder
(\cala)$ which means that the $Z(\cala)$-module $\gder(\cala)$ is the
$\cala$-dual of the bimodule
$\Omega^1_{\der}(\cala):\gder(\cala)=(\Omega^1_{\der}(\cala))^{\ast_{\cala}}$.
Thus $\os^1_{\der}(\cala)$ is the $\cala$-bidual bimodule of
$\Omega^1_{\der}(\cala)$, i.e. one has
$\os^1_{\der}(\cala)=(\Omega^1_{\der}(\cala))^{\ast_{\cala}\ast_{\cala}}$. This
is an obvious density result which implies in particular that
$\Omega^1_{\der}(\cala)$ is a diagonal bimodule; however this fact is obvious
since $\os^1_{\der}(\cala)$ is diagonal by definition $(\subset
\cala^{\der(\cala)})$.\\

Using
%% FOLLOWING LINE CANNOT BE BROKEN BEFORE 80 CHAR
%% FOLLOWING LINE CANNOT BE BROKEN BEFORE 80 CHAR
$\hom^{\cala}_{\cala}(\Omega^1_u(\cala),\cala)=\hom^{\cala}_{\cala}(\Omega^1_{\der}(\cala),\cala)$ one can characterize the pair $(\Omega^1_{\der}(\cala),d)$ consisting of the diagonal bimodule $\Omega^1_{\der}(\cala)$ and the derivation $d$ of $\cala$ into $\Omega^1_{\der}(\cala)$ by the following universal property \cite{mdv:pm2}: for any derivation $\delta$ of $\cala$ into a diagonal bimodule $M$ over $\cala$, there is a unique bimodule homomorphism $i_\delta:\Omega^1_{\der}(\cala)\rightarrow M$ such that $\delta=i_\delta\circ d$. This means that if $\delta$ is a derivation of $\cala$ into a diagonal bimodule $M$, the bimodule homomorphism $j_\delta:\Omega^1_u(\cala)\rightarrow M$ factorizes through the canonical bimodule homomorphism of $\Omega^1_u(\cala)$ onto $\Omega^1_{\der}(\cala)$. Recall that the underlying bimodule of the module of sections of a vector bundle over a manifold $V$ is diagonal and that $\Omega^1_{\der}(\calc^\infty(V))$ is the space of 1-forms on $V$!
, so the above result generalizes
a well known result of  differential geometry.\\

Let $X$ be a derivation of $\cala$, then one defines an antiderivation $i_X$ of
degree $-1$ of $\os_{\der}(\cala)$ by setting
$(i_X\omega)(X_1,\dots,X_{n-1})=\omega(X,X_1,\dots,X_n)$ for $\omega\in
\os^n_{\der}(\cala)$ and $X_i\in\gder(\cala)$. The mapping $X\mapsto i_X$ is an
operation, in the sense of H.~Cartan \cite{cart:01}, of the Lie algebra
$\gder(\cala)$ in the graded differential algebra $\os_{\der}(\cala)$, i.e. one
has $i_Xi_Y+i_Yi_X=0$ and, if one sets $L_X=\linebreak[4] di_X+i_Xd,
L_Xi_Y-i_YL_X=i_{[X,Y]}$ and $L_XL_Y-L_YL_X=L_{[X,Y]}$. Furthermore
$L_X\restriction\cala=X$ so $X\mapsto L_X$ is a Lie algebra homomorphism of
$\gder(\cala)$ into the derivations of degree zero of $\os_{\der}(\cala)$ which
extends the action of $\gder(\cala)$ on $\cala$. The differential subalgebra
$\Omega_{\der}(\cala)$ is stable by the $i_X,X\in\gder(\cala)$, so one has by
restriction an operation of $\gder(\cala)$ in $\Omega_{\der}(\cala)$. The
operation $X\mapsto i_X$ is of course the generalization of the interior
product (or contraction) of forms by vector fields while $L_X$ generalizes the
Lie derivative on forms.

\section{Noncommutative symplectic structures}

It is well known that the structural similarity between classical mechanics and
quantum mechanics is the most apparent if one uses the hamiltonian approach for
the former and that this is important for the problems of classical and
semiclassical limits. In this context the appropriate generalization of the
Poisson structures is also well known. {\sl A Poisson bracket on} $\cala$ is a
Lie algebra structure $(x,y)\mapsto \{x,y\}$ on $\cala$ satisfying
$\{x,yz\}=\{x,y\}z+y\{x,z\}$ for any elements $x,y$ and $z$ of $\cala$. Such a
Poisson bracket is {\sl real} if furthermore one has
$\{x,y\}^\ast=\{x^\ast,y^\ast\}$ for $x,y\in \cala$. For any $\cala$, there is
the standard real Poisson bracket $\{x,y\}=i[x,y]$ $(=i(xy-yx))$. Although this
bracket is trivial for a commutative algebra it is, up to a real factor, the
most common Poisson bracket occuring in quantum mechanics. In classical
hamiltonian mechanics, the Poisson bracket is associated with the symplectic
structure of the phase space. It is the aim of this section to describe the
generalization of symplectic structures for $\cala$ and to show its relevance
for quantum mechanics \cite{dv:2}, \cite{dvkm:1}, \cite{mad}.\\

The first thing to do is to generalize the notion of a nondegenerate two-form.
An element $\omega$ of $\os^2_{\der}(\cala)$ will be said to be {\sl
nondegenerate} if, for any $x\in\cala$, there is a derivation $\ham(x)\in
\gder(\cala)$ such that one has $\omega(X,\ham(x))=X(x)$ for any $X\in
\gder(\cala)$.  Notice that if $\omega$ is nondegenerate then $X\mapsto
i_X\omega$ is an injective linear mapping of $\gder(\cala)$ into
$\os^1_{\der}(\cala)$ but that the converse is not true; the condition for
$\omega$ to be nondegenerate is stronger than the injectivity of $X\mapsto
i_X\omega$. If $V$ is a manifold, an element $\omega\in
\os^2_{\der}(\calc^\infty(V))$ is an ordinary 2-form on $V$ and it is
nondegenerate in the above sense if and only if the 2-form $\omega$ is
nondegenerate in the classical sense (i.e. everywhere nondegenerate).\\

Let $\omega\in\os^2_{\der}(\cala)$ be nondegenerate, then for a given
$x\in\cala$ the derivation $\ham(x)$ is unique and $x\mapsto\ham(x)$ is a
linear mapping of $\cala$ into $\gder(\cala)$. Define then an antisymmetric
bilinear bracket on $\cala$ by $\{x,y\}=\omega(\ham(x),\ham(y))$. One has
$\{x,yz\}=\{x,y\}z+y\{x,y\}$ for $x,y,z\in \cala$,however the bracket
$(x,y)\mapsto\{x,y\}$ is a Lie bracket, (i.e. satisfies the Jacobi identity),
if and only if $d\omega=0$. A closed nondegenerate element $\omega$ of
$\os^2_{\der}(\cala)$ will be called {\sl a symplectic structure for} $\cala$.
Let $\omega$ be a symplectic structure for $\cala$, then the corresponding
bracket $(x,y)\mapsto \{x,y\}=\omega(\ham(x),\ham(y))$ is a Poisson bracket on
$\cala$ and one has $[\ham(x),\ham(y)]=\ham(\{x,y\})$, i.e. $\ham$ is a
Lie-algebra homomorphism of $(\cala,\{,\})$ into $\gder(\cala)$. If furthermore
$\omega$ is real, i.e. $\omega=\omega^\ast$, then this Poisson bracket is real
and $\ham(x^\ast)=(\ham(x))^\ast$ for any $x\in\cala$. We shall refer to the
above bracket as the Poisson bracket associated to the symplectic structure
$\omega$.\\

If $V$ is a manifold, a symplectic structure for $\calc^\infty(V)$ is just a
symplectic form on $V$. Since there are manifolds which do not admit symplectic
form, one cannot expect that an arbitrary $\cala$ admits a symplectic
structure.\\

Assume that $\cala$ has a trivial center $Z(\cala)=\Bbb C\bbbone$ and that all
its derivations are inner (i.e. of the form $ad(x),x\in\cala$). Then one
defines an element $\omega$ of $\os^2_{\der}(\cala)$ by setting
$\omega(ad(ix),ad(iy))=i[x,y]$. It is easily seen that $\omega$ is a real
symplectic structure for which one has $\ham(x)=ad(ix)$ and $\{x,y\}=i[x,y]$.
Although a little tautological, this construction is relevant for quantum
mechanics.\\

Let $\cala$ be, as above, a complex unital $\ast$-algebra with a trivial center
and only inner derivations and assume that there exists a linear form $\tau$ on
$\cala$ which is central, i.e. $\tau(xy)=\tau(yx)$, and normalized by
$\tau(\bbbone)=1$. Then one defines an element $\theta\in\os^1_{\der}(\cala)$
by $\theta(ad(ix))=x-\tau(x)\bbbone$. One has
$(d\theta)(ad(ix),ad(iy))=i[x,y]$, i.e. $\omega=d\theta$, so in this case the
symplectic form $\omega$ is exact. As examples of such algebras one can take
$\cala=M_n(\Bbb C)$, (a factor of type I$_n$), with $\tau=\frac{1}{n}$ trace,
or $\cala=\calr$, a von Neumann algebra which is a factor of type II$_1$ with
$\tau$  equal to the  normalized trace. The algebra $M_n(\Bbb C)$ is the
algebra of observables of a quantum spin $s=\frac{n-1}{2}$ while $\calr$ is the
 algebra used to describe the observables of an infinite assembly of quantum
spin; two typical types of quantum systems with no classical counterpart.\\

Let us now consider the C.C.R. algebra (canonical commutative relations)
$\cala_{CCR}$ \cite{dv:2}. This is the complex unital $\ast$-algebra generated
by two hermitian elements $q$ and $p$ satisfying the relation
$[q,p]=i\hbar\bbbone$. This algebra is the algebra of observables of the
quantum counterpart of a classical system with one degree of freedom. We keep
here the positive constant $\hbar$ (the Planck constant) in the formula for
comparison with classical mechanics, although the algebra for $\hbar\not=0$ is
isomorphic to the one with $\hbar=1$. We restrict here attention to one degree
of freedom to simplify the notations but the discussion extends easily to a
finite number of degrees of freedom. This algebra has again only inner
derivations and a trivial center so
$\omega(ad(\frac{i}{\hbar}x),ad(\frac{i}{\hbar}y))=\frac{i}{\hbar}[x,y]$
defines a symplectic structure for which $\ham(x)=ad(\frac{i}{\hbar}x)$ and
$\{x,y\}=\frac{i}{\hbar}[x,y]$ which is the standard quantum Poisson bracket.
In this case one can express $\omega$ in terms of the generators $q$ and $p$
and their differentials :
$$\omega=\sum_{n\geq 0}\left(\frac{1}{i\hbar}\right)^n \frac{1}{(n+1)!}[\dots
[dp,\underbrace{p],\dots,p]}_n[\dots[dq,\underbrace{q],\dots,q}_n]$$
Notice that this formula is meaningful because if one inserts  two derivations
$ad(ix),ad(iy)$ in it, only a finite number of terms contribute in the sum. For
$\hbar=0$, $q$ and $p$ commute and the algebra reduces to the algebra of
complex polynomial functions on the phase space $\Bbb R^2$. Furthermore the
limit of $\{x,y\}=\frac{i}{\hbar}[x,y]$ at $\hbar=0$ reduces to the usual
classical Poisson bracket as well known and, by using the above formula, one
sees that the formal limit of $\omega$ at $\hbar=0$ is $dpdq$.

\section{Derivations and Connections}

In this section $\calc$ is a complex unital commutative $\ast$-algebra and
$\cala$ is a complex unital $\ast$-algebra which is to be considered as the
noncommutative generalization of $\calc$. Our aim is to discuss the theory of
connections on the various objects which generalize the $\calc$-modules when
$\calc$ is replaced by $\cala$ in the framework of the differential calculus
based on derivations as generalization of vector fields \cite{mdv:pm1} (cf.
Section 0.4). In most parts of the following the involution is not involved and
therefore, in the definitions and results where the reality conditions do not
enter, one may assume that $\calc$ and its noncommutative counterpart $\cala$
are simply algebras (instead of $\ast$-algebras).\\

As generalizations of the category of $\calc$-modules when $\calc$ is replaced
by $\cala$ we consider the four following categories (cf. Section 0.2), the
category $\fc_{(0,0)}$ of $Z(\cala)$-modules, the category $\fc_{(1,0)}$ of
left $\cala$-modules, the category $\fc_{(0,1)}$ of right $\cala$-modules and
the category $\fc_{(1,1)}$ of central bimodules over $\cala$ i.e. of left
$\cala\otimes_{Z(\cala)} \cala^{op}$-modules. In each of these categories, one
has a direct sum and if $M$ is an object of any of these categories, it has a
canonical underlying structure of $Z(\cala)$-module. The labelling of these
categories by elements $\alpha=(i,j)$ of $\Bbb Z_2 \times \Bbb Z_2$ will be
very convenient to deal with the duality and tensor products. In $\Bbb Z_2
\times \Bbb Z_2$ one defines an involutive mapping $\alpha \mapsto \alpha'$ by
$(i,j)'=(1-i,1-j)$, i.e. $\alpha'=\alpha+(1,1)$. Correspondingly one has a
duality $M\mapsto M'$ of $\fc_\alpha$ into $\fc_{\alpha'}$, where $M'=M^\ast$
if $M$ is a left or right $\cala$-module and $M'=M^{\ast_{\cala}}$ if $M$ is a
$Z(\cala)$-module or a central bimodule over $\cala$. Another bit of notation
will be convenient; we set $A_0=Z(\cala)$ and $A_1=\cala$. Using this notation,
an object of $\fc_{(i,j)}$ is a $(A_i,A_j)$-bimodule (of a specific kind) and
we can define tensor products $\fc_{(i,j)} \times \fc_{(j,k)} \rightarrow
\fc_{(i,j)} \tilde{\otimes} \fc_{(j,k)} \subset \fc_{(i,k)}$ by
$M\tilde{\otimes} N=M\otimes_{A_j}N$ if $M$ is an object of $\fc_{(i,j)}$ and
$N$ an object of $\fc_{(j,k)}$ (one verifies that $M\tilde{\otimes}N$ is then
an object of $\fc_{(i,k)}$).\\

Let $M$ be an object of $\fc_{(i,j)}$. A {\sl connection on} $M$ is a linear
mapping $\nabla$, $X\mapsto \nabla_X$, of $\gder(\cala)$ into the linear
endomorphism of $M$ such that one has for any $m\in M$ and any
$X\in\gder(\cala)$
\[
\left\{
\begin{array}{l}
\nabla_{zX}(m)=z\nabla_X(m),\ \ \forall z\in Z(\cala)\\
\nabla_{X}(a_ima_j)=X(a_i)ma_j+a_i\nabla_X(m)a_j+a_imX(a_j),\ \ \forall a_i\in
A_i, \forall a_j\in A_j
\end{array}
\right.
\]
remembering that $M$ is canonically a $Z(\cala)$-module and that since
$Z(\cala)=A_0$ is stable by $\gder(\cala),\gder(\cala)$ acts by derivations on
$Z(\cala)=A_0$ and on $\cala=A_1$. It should be stressed that elements of
$A_0=Z(\cala)$ can be moved to the other side. Given $\nabla$ as above, {\sl
the curvature} $R$ {\sl of} $\nabla$ is the bilinear antisymmetric mapping
$(X,Y)\mapsto R_{X,Y}$ of $\gder(\cala)\times \gder(\cala)$ into the linear
endomorphisms of $M$ defined by
$R_{X,Y}(m)=\nabla_X(\nabla_Y(m))-\nabla_Y(\nabla_X(m))-\nabla_{[X,Y]}(m),\
\forall X,Y\in \gder(\cala),\ \forall m\in M$. One has
$R_{zX,Y}(m)=zR_{X,Y}(m)$ and $R_{X,Y}(a_ima_j)=a_iR_{X,Y}(m)a_j,\ \forall m\in
M$, $\forall X,Y\in \gder(\cala)$, $\forall z\in Z(\cala)$, $\forall a_i\in
A_i$, $\forall a_j\in A_j$. More precisely, $R$ is an antisymmetry
$Z(\cala)$-bilinear mapping of $\gder(\cala)\times \gder(\cala)$ into the
$Z(\cala)$-module $\hom_{\fc_{(i,j)}}(M,M),(\hom_{\fc_{(i,j)}}$ being the
morphisms in $\fc_{(i,j)}$).\\

There is an obvious connection $\nabla_1\oplus \nabla_2$ on the direct sum
$M_1\oplus M_2$ of two objects $M_1$ and $M_2$ of $\fc_{(i,j)}$ equipped with
connections $\nabla_1$ and $\nabla_2$.\\

Let $M$ be an object of $\fc_{(i,j)}$ then its dual $M'$ is an element of
$\fc_{(i,j)'}$ and we denote by $(m,m')\mapsto <m,m'>\in \cala$ the bilinear
duality bracket obtained by evaluation, $<,> : M\times M'\rightarrow \cala$.
Then, for any connection $\nabla$ on $M$, there is unique {\sl dual connection}
$\nabla'$ on $M'$ such that  $X(<m,m'>)=<\nabla_X(m),m'>+<m,\nabla'_X(m')>$,
$\forall m\in M, \forall m'\in M'$ and $\forall X\in \gder(\cala)$. Indeed the
above equality defines $\nabla'$ uniquely and one checks that it is a
connection. In general, the mapping $\nabla\mapsto\nabla'$ is not injective nor
surjective. However if the canonical mapping of $M$ into its bidual $M''$,
(which is a morphism of $\fc_{(i,j)}$), is injective, then $\nabla''$ is an
extension of $\nabla$ and therefore $\nabla\mapsto \nabla'$ is injective and of
course bijective whenever $M=M''$. An object $M$ of $\fc_{(ij)}$ will be called
{\sl diagonal} if the canonical morphism of $M$ in $M''$ is injective. This
generalizes the notion introduced in Section 0.2 (for $\fc_{(1,1)}$) and the
terminology is suggested by the following. The algebra $\cala$ itself can be
considered as an object of $\fc_{(i,j)}$ when it is equipped with the canonical
corresponding underlying structure and the same is true for $\cala^I$ where $I$
is an arbitrary set, since $\fc_{(i,j)}$ has arbitrary products and, more
generally, arbitrary projective limits. Then $M$ is diagonal if and only if
there is an injective $\fc_{(i,j)}$-morphism of $M$ into $\cala^I$, for some
set $I$. Finally let us notice that any projective limit of diagonal objects is
diagonal and that dual objects are diagonal, i.e. if $M$ is an object of
$\fc_{(i,j)}$, then its dual $M'$ is a diagonal object of $\fc_{(i,j)'}$.

Let $M_1$ be an object of $\fc_{(i,j)}$ and $M_2$ be an object of $\fc_{(j,k)}$
and let $\nabla^1$ be a connection on $M_1$ and $\nabla^2$ be a connection on
$M_2$. Then, for any $X\in\gder(\cala)$, $D_X=\nabla^1_X\otimes
id_{M_{2}}+id_{M_{1}}\otimes \nabla^2_X$ is such that that it maps into itself
the subspace of $M_1\otimes M_2$ generated by the elements $m_1a_j\otimes
m_2-m_1\otimes a_jm_2$, with $m_1
\in M_1$, $m_2\in M_2$ and $a_j\in A_j$. It follows that the $D_X$ pass to the
quotient and define linear endomorphisms $\nabla_X$ of $M_1\tilde{\otimes}M_2$
and one verifies that $\nabla$ so defined is a connection on the object
$M_1\tilde{\otimes}M_2$ of $\fc_{(i,k)}$. This connection will be refered to as
{\sl the tensor product of $\nabla_1$ and $\nabla_2$}.

Thus we have defined connections on $Z(\cala)$-modules, on left and right
$\cala$-modules and on central bimodules over $\cala$ and we have also defined
dual and tensor product of such connections. Let us now come to the problems of
reality for such connections. As pointed out in Section 0.2, the notion of
reality makes sense only for $\ast$-modules over  $Z(\cala)$ or for
$\ast$-bimodules over $\cala$. So let $M$ be either a $\ast$-module over
$Z(\cala)$ or a $\ast$-bimodule over $\cala$ which is central. If $\nabla$ is a
connection on $M$ one can define another one $\nabla^\ast$, {\sl its
conjugate}, by setting $\nabla^\ast_X(m)=(\nabla_{X^{\ast}}(m^\ast))^\ast$ and
$\nabla$ will be said to be a {\sl real} connection if $\nabla=\nabla^\ast$.
Let $M'$ be the dual of $M$, i.e. $M'=M^{\ast_{\cala}}$ in this case, then
there is a unique involution $m'\mapsto m'{^\ast}$ on $M'$ such that
$<m,m'>{^\ast}=\linebreak[4] <m^{\ast},m'{^\ast}>$ and, equipped with this
involution, $M'$ is a (central) $\ast$-bimodule over $\cala$ if $M$ is a
$\ast$-module over $Z(\cala)$ or a $\ast$-module over $Z(\cala)$ if $M$ is a
central $\ast$-bimodule over $\cala$. Furthermore, one has
$(\nabla^\ast)'=(\nabla')^\ast$, so the dual connection of a real connection is
real.

\section{Linear connections}

In classical differential geometry, a connection on the tangent bundle, or
equivalently on the cotangent bundle, of a manifold is usually called a linear
connection. Although this terminology is a little misleading, we shall
nevertheless use it for the corresponding noncommutative generalizations.
Within the framework of Section~0.4 and Section~0.6, one sees that there are
three natural definitions of such generalizations. First a connection on
$\Omega^1_{\der}(\cala)$, second a connection on $\gder(\cala)$ and third a
connection on $\os^1_{\der}(\cala)$. However, as explained in Section 0.4,
$\Omega^1_{\der}(\cala)$ is a diagonal bimodule with $\gder(\cala)$ as
$\cala$-dual, i.e. $\gder(\cala)=(\Omega^1_{\der}(\cala))'$ with the notation
of Section 0.6, and $\os^1_\der(\cala)$ is the $\cala$-dual of $\gder(\cala)$,
i.e. $\os^1_{\der}(\cala)=(\gder(\cala))'=(\Omega^1_{\der}(\cala))''$.
Therefore, it follows from the discussion of the previous section that, by
duality, there is an injective mapping of the (affine) space of connections on
$\Omega^1_{\der}(\cala)$ into the space of connections on $\gder(\cala)$ and
that there is also on injective mapping of the space of connections on
$\gder(\cala)$ into the space of connections on $\os^1_{\der}(\cala)$. Thus all
these connections may be imbedded into the connections on
$\os^1_{\der}(\cala)$. A real connection on $\os^1_{\der}(\cala)$ will be
called {\sl a linear connection on} $\cala$. The connections on $\gder(\cala)$
form a subclass of connections on $\os^1_{\der}(\cala)$ and an even smaller
subclass consists of connections on $\Omega^1_{\der}(\cala)$. Given a
connection $\nabla$ on $\os_{\der}(\cala)$, one defines a bimodule homomorphism
$T:\os^1_{\der}(\cala)\rightarrow \os^2_{\der}(\cala)$, {\sl its torsion}, by
setting $(T\omega)(X,Y)=(d\omega)(X,Y)-\nabla_X(\omega)(Y)+\nabla_Y(\omega)(X)$
for $X,Y\in \gder(\cala)$ and $\omega\in\os^1_{\der}(\cala)$. If $\nabla$ comes
from a connection on $\Omega^1_{\der}(\cala)$, (by biduality), $T$ restricted
to $\Omega^1_{\der}(\cala)$, is a bimodule homomorphism of
$\Omega^1_{\der}(\cala)$ into $\Omega^2_{\der}(\cala)$. If $\nabla$ is the dual
of a connection, again denoted by $\nabla$, on $\gder(\cala)$, its torsion can
be identified with the $Z(\cala)$-bilinear antisymmetric mapping $T$ of
$\gder(\cala)\times \gder(\cala)$ into $\gder(\cala)$ defined by
$T(X,Y)=\nabla_X(Y)-\nabla_Y(X)-[X,Y]$, $\forall X,Y\in \gder(\cala)$. For a
more complete discussion as well as for the notion of Levi-Civita connection of
a generalization of pseudo-riemannian metric, we refer to \cite{mdv:pm1}.

\section{Conclusion: General differential calculi}

The above notions of connections are natural ones when one uses the
differential calculus based on derivations as generalization of vector fields.
However, for some purposes, (see e.g. in \cite{connes:03}), it is useful to use
other differential calculi and therefore, it is natural to ask for a definition
of connections adapted to such calculi. Let $\Omega$ be a differential calculus
over $\cala$. There is then a well known useful definition of an
$\Omega$-connection on  a left (or right) $\cala$-module \cite{connes:01}. The
problem arises when one tries to define an $\Omega$-connection on a bimodule
over $\cala$ such as $\Omega^1$. This problem is unavoidable if one wishes to
generalize linear connections since the natural structure of $\Omega^1$ is that
of a  bimodule. Some authors, e.g. \cite{cfg}, define a connection on
$\Omega^1$ to be a left module $\Omega$-connection on $\Omega^1$. Besides the
fact that it is  unnatural to privilege part of a bimodule structure, this
definition has two drawbacks if one thinks of it as a generalization of linear
connections. First, one cannot introduce then the notion of reality which
generalizes the classical notion of reality of a linear connection in
differential geometry because the involution of $\Omega^1$ is linked to its
bimodule structure, (so the conjugate of a left $\cala$-module connection on
$\Omega^1$ is rather a right $\cala$-module connection). Second, one cannot, in
general, define the tensor product over $\cala$ of such a connection with a
connection say on a left $\cala$-module although it is very desirable to have
such a tensor product, e.g. for the description of the generalization of the
classical coupling of gravitation with a field coupled to a Yang-Mills field. A
definition of linear  connections for general differential calculi which takes
into account the complete bimodule structure of $\Omega^1$ has been proposed by
J. Mourad \cite{jm} and further generalized to other bimodules \cite{mdv:m}.
This definition involves a generalization of the permutations in tensor
products and with it the question of reality can be addressed. Furthermore, in
this framework, tensor products of connections are defined straightforwardly
and it has been recently shown \cite{bmds}, (see in appendix A of \cite{bmds}),
that conversely, in order that tensor products of connections exist in a very
general sense, one has to use this definition of connections for bimodules.

\bibliographystyle{amsplain}

\end{document}